\documentclass[11pt]{article}

\usepackage{amsmath,amssymb,amstext,amsthm,graphicx}
\usepackage[sort&compress,round,authoryear]{natbib}
\usepackage[left=3.15cm,right=3.15cm,top=2.9cm,bottom=2.9cm]{geometry}

\newtheorem{theorem}{Theorem}
\newtheorem{corollary}[theorem]{Corollary}
\newtheorem{proposition}[theorem]{Proposition}

\newtheorem{definition}[theorem]{Definition}
\newcounter{remcount}

\newcounter{example}

\newcommand{\ind}{\mbox{$\perp \kern-5.5pt \perp$}}

\numberwithin{equation}{section}

\begin{document}

\title{{\Large\bf Multiple Solutions to the Likelihood Equations in
   the Behrens-Fisher Problem}} 
\author{\normalsize 
  Mathias Drton\footnote{This material is based upon work supported by the National Science
Foundation under Grant No.~0505612.}\\
  \em\normalsize
  Department of Statistics, University of Chicago,\\
  \em\normalsize 5734 S. University
  Avenue, Chicago, Illinois 60637, U.S.A.\\
  \rm\normalsize Email: {\tt drton@galton.uchicago.edu}
} 

\date{} 
\maketitle

\renewcommand{\labelenumi}{(\roman{enumi})}
\renewcommand{\thefootnote}{\fnsymbol{footnote}}
\renewcommand{\thetable}{\arabic{table}}
\renewcommand{\labelenumi}{(\roman{enumi})}

\maketitle

\begin{abstract} \normalsize
  The Behrens-Fisher problem concerns testing the equality of the means of
  two normal populations with possibly different variances.  The null
  hypothesis in this problem induces a statistical model for which the
  likelihood function may have more than one local maximum.  We show that
  such multimodality contradicts the null hypothesis in the sense that if
  this hypothesis is true then the probability of multimodality converges
  to zero when both sample sizes tend to infinity.
  Additional results include a finite-sample bound on the probability of
  multimodality under the null and asymptotics for the probability of
  multimodality under the alternative.\\

\end{abstract}

\section{Introduction}
\label{sec:intro}

The Behrens-Fisher problem is concerned with testing
\[
H_0:\mu_X=\mu_Y \quad\text{vs.}\quad H_1:\mu_X\not=\mu_Y,
\]
where $\mu_X$ and $\mu_Y$ are the means of two normal populations with
possibly different variances $\sigma_X^2$ and $\sigma_Y^2$.  An interesting
aspect of the problem is that the likelihood equations for the model
induced by $H_0$ may have more than one solution.  In fact, with
probability one, there will be either one or three solutions with the two
cases corresponding to one or two local maxima of the likelihood function.
According to simulations of \cite{sugiura}, three solutions to the
likelihood equations occur infrequently if the observations are drawn from
a distribution in $H_0$.  In this note we provide an explanation for this
rare occurrence of multiple solutions by proving that under $H_0$ the
probability of this event converges to zero when $n$ and $m$ tend to
infinity (Corollary~\ref{cor:asy}).  This and more general large-sample
results about the probability of multiple solutions
(Proposition~\ref{prop:disc-pos-neg} and Theorem~\ref{thm:critical-smooth})
are based on two observations.  First, solving the likelihood equations
amounts to solving one cubic polynomial equation.  Second, the number of
real roots of a cubic can be determined using the cubic discriminant.  The
discriminant criterion also allows us to derive a finite-sample bound on
the null probability of multiple solutions to the likelihood equations
(Proposition~\ref{prop:t-distribution-bound}).

While arguments can be given for using a likelihood ratio-based test
instead of Welch's approximate t-test in the Behrens-Fisher problem
\citep{jensen:1992}, the latter test is widely used in practice and avoids
maximization of the likelihood function under the null hypothesis.  In that
sense the practical implications of our study are perhaps not immediate.
However, in more general models involving heteroscedastic structures
statistical practice often relies on likelihood ratio tests that do require
solving the maximization problem.  Our results provide geometric intuition
about this problem in the simple univariate Behrens-Fisher model
(Figure~\ref{fig:curves}), for which it holds, rather reassuringly, that
the likelihood function for the null model is asymptotically unimodal if
the model is correctly specified.  It would be interesting to obtain
generalizations of this fact for other, more complicated models.

\section{Solving the likelihood equations}
\label{sec:lik-equation}

We start out by deriving a convenient form of the likelihood equations for
the three-parameter model induced by the null hypothesis in the
Behrens-Fisher problem.

\subsection{A cubic equation}
\label{subsec:cubic}


Let $X_1,\dots,X_{n} \sim N(\mu_X,\sigma_X^2)$ and $Y_1,\dots,Y_{m} \sim
N(\mu_Y,\sigma_Y^2)$ be two independent normal samples.  Under the null
hypothesis $H_0$, $\mu_X$ is equal to $\mu_Y$ and we denote this common
mean by $\mu$.  The log-likelihood function for the null model can be
written as
\begin{align*}
  \ell(\mu,\sigma_X^2,\sigma_Y^2) =&  -\frac{n+m}{2} \log (2\pi)-
  \frac{n}{2}\log(\sigma_X^2) 
  -\frac{m}{2}\log(\sigma_Y^2)\\  &
  -\frac{n}{2} \left[ \frac{\hat\sigma_X^2+(\bar X-\mu)^2}{\sigma_X^2}\right]
  -\frac{m}{2} \left[ \frac{\hat\sigma_Y^2+(\bar Y-\mu)^2}{\sigma_Y^2}\right].
\end{align*}
Here, $\bar X$ and $\bar Y$ are the two sample means, and
\[
\hat\sigma_X^2=\frac{1}{n}\sum_{i=1}^n(X_i-\bar X)^2
\]
is the empirical variance for the first sample; the second empirical
variance $\hat\sigma_Y^2$ is defined analogously.  If $\min(n,m)\ge 2$,
then both $\hat\sigma_X^2$ and $\hat\sigma_Y^2$ are positive with
probability one.  This sample size condition will be assumed throughout.

The partial derivatives of the log-likelihood function are
\begin{align*}
  \frac{\partial\ell}{\partial\mu}&=
  \frac{n(\bar X-\mu)}{\sigma_X^2}+
  \frac{m(\bar Y-\mu)}{\sigma_Y^2}\\ 
\intertext{and}
  \frac{\partial\ell}{\partial\sigma_X^2}&=
  -\frac{n}{2\sigma_X^2}+\frac{n[\hat\sigma_X^2+(\bar X-\mu)^2]}{2\sigma_X^4};
\end{align*}
the partial derivative for $\sigma_Y^2$ is analogous.  Let $r_n=n/m$.  Then
the likelihood equations obtained by setting the three partial derivatives
to zero are equivalent to the polynomial equations
\begin{gather}
  \label{eq:1}  
  r_n(\bar X-\mu)\sigma_Y^2+(\bar Y-\mu)\sigma_X^2 = 0,\\
  \label{eq:2}
  \sigma_X^2 = (\bar X-\mu)^2 + \hat\sigma_X^2, \\
  \label{eq:3}
  \sigma_Y^2 = (\bar Y-\mu)^2 + \hat\sigma_Y^2.
\end{gather}
Here, equivalence means that the two solution sets are almost surely equal.
Plugging the expressions for $\sigma_X^2$ and $\sigma_Y^2$ from
(\ref{eq:2}) and (\ref{eq:3}) into (\ref{eq:1}) yields the cubic equation
\begin{equation}
  \label{eq:cubic}
  f(\mu)=a_3\mu^3+a_2\mu_2+a_1\mu+a_0=0
\end{equation}
with
\begin{align*}
  \nonumber
  a_3 &= 1+r_n,\\
  a_2 &= - (2\bar X+\bar Y) - r_n(2\bar Y+\bar X),\\
  a_1 &= \bar X^2  + 2(1+r_n)\bar X\bar Y +
  r_n\bar Y^2 +\hat\sigma_X^2  + r_n\hat\sigma_Y^2, \quad\text{and} \\
  a_0 &=- \bar X^2\bar Y - r_n\bar Y^2\bar X
  -\hat\sigma_X^2\bar Y  - r_n\hat\sigma_Y^2 \bar X.
\end{align*}
Hence, the maximum likelihood estimator $\hat\mu$ can be computed in closed
form by solving the univariate cubic equation (\ref{eq:cubic}).
We remark that the manipulations leading to the polynomial
equations (\ref{eq:2}), (\ref{eq:3}) and (\ref{eq:cubic}) form a trivial
case of a computation of a {\em lexicographic Gr\"obner basis}
\citep[p.~86]{ascb}.  

\subsection{The discriminant}
\label{subsec:disc}

A quadratic polynomial $a_2x^2+a_1x+a_0$ in the indeterminate $x$ may have
no, one, or two (distinct) real roots.  Which one of the three cases
applies is determined by the sign of the discriminant $a_1^2-4a_0a_2$.  In
the Behrens-Fisher problem we are led to the cubic polynomial $f$ in
(\ref{eq:cubic}).  A cubic always has at least one real root, and so we
would like to know whether it has one, two, or three real roots.  This can
again be decided based on the sign of the {\em discriminant\/}, which for
the cubic takes the form
\[
\Delta = a_1^2a_2^2-4a_0a_2^3-4a_1^3a_3+18a_0a_1a_2a_3-27a_0^2a_3^2.
\]
If $\Delta>0$ then $f$ has three distinct real roots.  If $\Delta<0$ then
$f$ has a unique real root and two complex ones.  If $\Delta=0$, then $f$
may have one real root of multiplicity three or two distinct real roots of
which one has multiplicity two.  These and more general results on
discriminants can be found for example in \citet[\S4.1-2]{basu:2003}.

The coefficients $a_0$, $a_1$ and $a_2$ of $f$ in (\ref{eq:cubic}) are
random variables with a continuous distribution and $a_3$ is a constant.
Consequently, $\Delta$ is also a continuous random variable such that the
event $\{\Delta=0\}$ occurs with probability zero.  In other words the
Behrens-Fisher likelihood equations almost surely have one or three real
solutions.

The discriminant $\Delta$ is a homogeneous polynomial of degree 6 in $\bar
X$, $\bar Y$, $\hat\sigma_X$ and $\hat\sigma_Y$, and depends on $\bar X$
and $\bar Y$ only through their difference.
However, with probability
one, the sign of $\Delta$ depends only on $r_n$ and the two ratios
$\hat\gamma=\hat\sigma_X/\hat\sigma_Y$ and $\hat\delta=(\bar X-\bar Y)/\hat\sigma_Y$.
This follows because $\Delta = \hat\sigma_Y^6\cdot D$
with
\begin{equation}
  \label{eq:D}
  \begin{split}
    D&=\hat\delta^6r_n^2-2\hat\delta^4\bigg[\hat\gamma^2(2+2r_n-r_n^2)
    +(2r_n^3+2r_n^4-r_n^2)\bigg] \\
    &- \hat\delta^2
    \bigg[\hat\gamma^4(8+8r_n-r_n^2)+(8r_n^4+8r_n^3-r_n^2)
    -2\hat\gamma^2(10r_n+19r_n^2+10r_n^3)
    \bigg]\\
    &-4(1+r_n)(r_n + \hat\gamma^2)^3.
  \end{split}
\end{equation}
While $\Delta$ (and $D$) remain unchanged if $n$ and $m$ are replaced by
$\bar n$ and $\bar m$ with $r_n=n/m=\bar n/\bar m=r_{\bar n}$, such a
change of sample sizes affects the sampling distribution of $\Delta$ (and
$D$) and thus the probability of multiple solutions to the likelihood
equations.  We remark that instead of working with $\hat\delta$ one could
work with the more symmetric quantities
\[
\frac{\bar X-\bar Y}{\sqrt{\hat\sigma_X^2+\hat\sigma_Y^2}}=
\hat\delta \cdot \frac{1}{\sqrt{\hat\gamma^2+1}}
\]
or
\[
\frac{\bar X-\bar Y}{\sqrt{\hat\sigma_X^2+r_n\hat\sigma_Y^2}}=
\hat\delta \cdot \frac{1}{\sqrt{\hat\gamma^2+r_n}}.
\]
However, such a substitution would lead to an increased degree in the
analog of (\ref{eq:D}) such that we keep working with $\hat\delta$ in the
sequel.

\begin{figure}[t]
  \centering
  \begin{tabular}{cc}
  (a) \includegraphics[width=6cm]{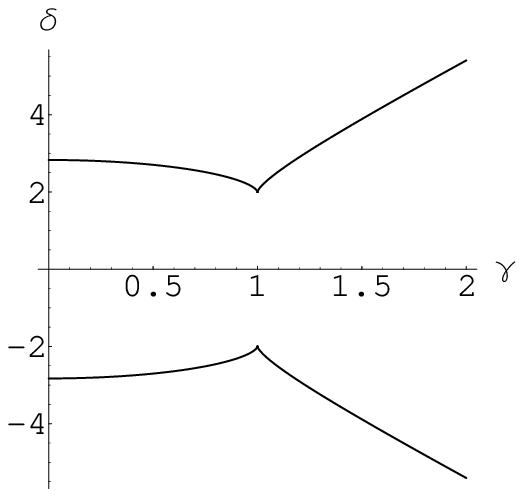}
  &$\quad$
  (b) \includegraphics[width=6cm]{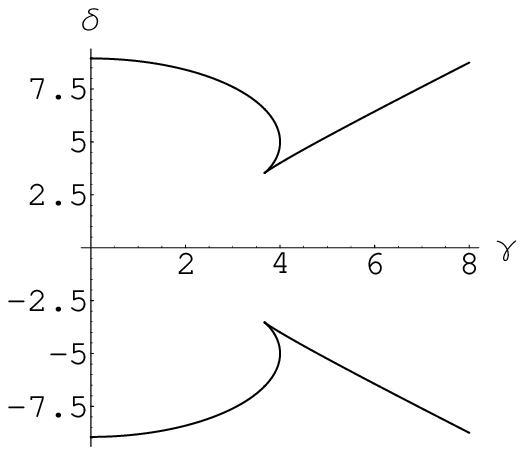}
  \end{tabular}
  \caption{Algebraic curve defined by the polynomial $D_{r_n}$ that is
    derived from the discriminant: (a) $r_n=1$ and (b) $r_n=4$.  In each
    plot, points $(\gamma,\delta)$ between the two curves correspond to a
    unique real root to the Behrens-Fisher likelihood equations. Points
    above and below the curves correspond to three distinct real roots.}
  \label{fig:curves}
\end{figure}

For any given value of $r_n$, the polynomial $D=D_{r_n}$ in the
indeterminates $\hat\gamma$ and $\hat\delta$ defines an algebraic curve.
Figure \ref{fig:curves} shows two examples of these curves over the
statistically relevant region with $\hat\gamma> 0$.  By symmetry, the curve
for $r_n=1$ has four cusp points at $(\hat\gamma,\hat\delta)=(\pm 1,\pm
2)$; the cusps for $r_n=4$ are at
$(\hat\gamma,\hat\delta)=(\pm\sqrt{27/2},\pm\sqrt{25/2})$.  In general, the
four cusps have coordinates
\begin{equation}
  \begin{split}
    \label{eq:cusp-coordinates}
    \hat\gamma&=\pm\frac{(2r_n+1)\sqrt{(r_n+2)r_n(2r_n+1)}}{(r_n+2)^2},\\
    \hat\delta&=\pm\frac{3(1+r_n)\sqrt{3(r_n+2)r_n}}{(r_n+2)^2}.
  \end{split}
\end{equation}
The curve has two asymptotes, namely, $r_n\hat\delta=\pm 2\sqrt{1+r_n}\,
\hat\gamma$.

The two respective branches of the curves in Figure \ref{fig:curves}
enclose the region $\{D<0\}$, which contains the (horizontal)
$\hat\gamma$-axis.  Over this region the discriminant $\Delta$ is negative
and the Behrens-Fisher likelihood equations have a unique real root.
Clearly, neither the region $\{D<0\}$ nor the region $\{D>0\}$ need to be
convex.  
When fixing $\hat\gamma$ and increasing $\hat\delta$ then $D$ will
eventually remain positive because the leading term of $D$, when viewed as
a univariate polynomial in $\hat\delta$, is $r_n^2\hat\delta^6$ with
$r_n^2>0$.  This means that bimodal likelihood functions arise when the
difference between the means $\bar X$ and $\bar Y$ of the two samples is
large compared to the empirical variances $\hat\sigma_X^2$ and
$\hat\sigma_Y^2$.  However, as can be seen in Figure \ref{fig:curves}(b)
with $r_n=4$, there may exist values of $\hat\gamma$ such that the values
of $\hat\delta$ corresponding to unimodal likelihood functions do not form
an interval.

\subsection{Finite-sample bound}
\label{sec:bound}

A finite-sample study of the probability of one versus three solutions to
the Behrens-Fisher likelihood equations seems difficult.  However, under
the null hypothesis, we can give a very simple bound.

\begin{proposition}
  \label{prop:t-distribution-bound}
  Let the random variable $T$ have a t-distribution with $m-1$ degrees of
  freedom.  Let $\gamma=\sigma_X/\sigma_Y$.  If the null hypothesis $H_0$
  is true, i.e., if $\mu_X=\mu_Y$, then the probability of three distinct
  real solutions to the Behrens-Fisher likelihood equations is smaller than
  \[
  P\left( |T|> \sqrt{m-1}\cdot \frac{3(1+r_n)r_n\sqrt{3(r_n+2)}}{(r_n+2)^2\sqrt{
        \gamma^2+r_n}} \right).
  \]
\end{proposition}
\begin{proof}
  Three solutions occur if $(\hat\gamma,\hat\delta)$ falls in the region
  $\{D>0\}$.  This region is strictly contained in the region of pairs
  $(\hat\gamma,\hat\delta)$ that have $|\hat\delta|>c_n$ with
  $$c_n=\frac{3(1+r_n)\sqrt{3(r_n+2)r_n}}{(r_n+2)^2};$$
  compare Figure
  \ref{fig:curves} and (\ref{eq:cusp-coordinates}).  Hence, $P(D>0)$ is
  smaller than $P(|\hat\delta|>c_n)$.  Under $H_0$,
  \[
  \frac{\sigma_Y}{\sqrt{\frac{m}{m-1}}\cdot
    \sqrt{\frac{\sigma_X^2}{n}+\frac{\sigma^2_Y}{m}}} 
  \cdot \hat\delta 
  \]
  is distributed according to the t-distribution with $m-1$ degrees of
  freedom.  Expressing the event $\{|\hat\delta|>c_n\}$ in terms of this
  t-random variable yields the claim.
\end{proof}

Suppose, for example, that the samples are of equal size with the standard
deviation $\sigma_X$ being half the standard deviation $\sigma_Y$, i.e.,
$r_n=1$ and $\gamma=1/2$.  Then, by Proposition
\ref{prop:t-distribution-bound}, the probability of three distinct real
solutions to the Behrens-Fisher likelihood equations is
smaller than 
0.023 if $n=m=5$, 0.00045 if $n=m=10$ and 0.00001 if $n=m=15$.  Hence,
despite its crude nature, the bound informs us that the probabilities are
small.  Monte Carlo simulations suggest that the three considered
probabilities are in fact a factor 10 or more smaller than the stated
bounds.

\section{Large-sample results}
\label{sec:large-samples}

We begin our study of the large-sample behaviour of the 
likelihood equations with the
case when the discriminant converges almost surely to a non-zero limit.



\begin{proposition}
  \label{prop:disc-pos-neg}
  Suppose $\min(n,m)\to\infty$ such that $r_n=n/m\to r\in(0,\infty)$.  Let
  $\delta=(\mu_X-\mu_Y)/\sigma_Y$ and $\gamma=\sigma_X/\sigma_Y$.  Define
  $D_r(\gamma,\delta)$ to be the quantity obtained from $D$ in (\ref{eq:D})
  by replacing $r_n$ by $r$ and $(\hat\gamma,\hat\delta)$ by
  $(\gamma,\delta)$.
  \begin{enumerate}
  \item[(i)] If $D_r(\gamma,\delta)<0$, then the probability that the
    Behrens-Fisher likelihood equations have exactly one real solution
    converges to one.
  \item[(ii)] If $D_r(\gamma,\delta)>0$, then the probability that the
    Behrens-Fisher likelihood equations have three distinct real solutions
    converges to one.
  \end{enumerate}
\end{proposition}
\begin{proof}
  The polynomial $D$ is a continuous function of 
  $\bar X$, $\bar Y$, $\hat\sigma_X^2$ and $\hat\sigma_Y^2$.  Applying laws
  of large numbers to the four random variables, we find that
  $D_{r_n}(\hat\gamma,\hat\delta)$ converges almost surely to
  $D_r(\gamma,\delta)$.  In case (i), $D_r(\gamma,\delta)$ is negative, and
  thus $P(D_{r_n}(\hat\gamma,\hat\delta)<0)$ converges to one, which implies the
  claim.  Case (ii) is analogous. 
\end{proof}

The next result is a corollary to both Propositions
\ref{prop:t-distribution-bound} and \ref{prop:disc-pos-neg}.


\begin{corollary}
  \label{cor:asy}
  Suppose $H_0$ is true, i.e., $\mu_X=\mu_Y=\mu$.  If $\min(n,m)\to\infty$
  and $r_n=n/m\to r\in(0,\infty)$, then the probability that the
  Behrens-Fisher likelihood equations have exactly one real solution
  converges to one.
\end{corollary}
\begin{proof}
  If $\mu_X=\mu_Y$, then $\delta=0$ and the claim follows from Proposition
  \ref{prop:disc-pos-neg} because $D_r(\gamma,0) = -4(1+r)(r + \gamma^2)^3$
  is negative.
\end{proof}

Proposition~\ref{prop:disc-pos-neg} does not apply to the situation when $D_r$
is zero.  However, these critical cases can be studied using asymptotics
similar to those encountered with likelihood ratio tests.  
The resulting asymptotic probabilities will depend on whether or not the
point $(\gamma,\delta)$ forms a singular point of the curve defined by the
vanishing of $D_r$.

\begin{definition}
  Let $h$ be a polynomial in the ring 
  of polynomials
  in the indeterminates $x_1$ and $x_2$ with real coefficients.  Let $V(h)$
  be the algebraic curve $\{x\in\mathbb{R}^2\mid h(x)=0\}$.  A point $x\in
  V(h)$ is a singular point if the gradient $\nabla h(x)$ is zero.
\end{definition}

Our curve of interest, $V(D_r)$, has four singular points whose coordinates
were given in (\ref{eq:cusp-coordinates}); recall Figure \ref{fig:curves}.
All other points on $V(D_r)$ are non-singular.

We will show that the critical behaviour of the number of real roots to the
Behrens-Fisher likelihood equations is determined by the local
geometry of the curve $V(D_r)$ at the true parameter values
$(\gamma,\delta)$.  This geometry is captured in the tangent cone.

\begin{definition} 
  The tangent cone of $V\subseteq\mathbb{R}^2$ at $x\in \mathbb{R}^2$ is
  the set of vectors that are limits of sequences $\alpha_n(x_n-x)$, where
  $\alpha_n$ are positive reals and $x_n\in V$ converge to $x$.
\end{definition}

The tangent cone, which is a closed set, is indeed a {\em cone\/}.  This
means that if $\tau$ is in the tangent cone then so is the half-ray
$\{\lambda\tau \mid \lambda\ge0\}$.  

\begin{theorem}
  \label{thm:critical-smooth}
  Suppose that $\min(n,m)\to\infty$ and $r_n=r+o(1/\sqrt{n})$.  Let
  $\gamma>0$.
  \begin{enumerate}
  \item[(i)] If $(\gamma,\delta)$ is a non-singular point of the curve
    $V(D_r)$, then the probability of exactly one real solution as well as
    the probability of three distinct real solutions to the Behrens-Fisher
    likelihood equations converge to $1/2$.
  \item[(ii)] If $(\gamma,\delta)$ is one of two singular points of the
    curve $V(D_r)$ that have $\gamma>0$, then the probability of exactly
    one real solution converges to one.
  \end{enumerate}
\end{theorem}
\begin{proof}
  We first show that the asymptotic probability can be obtained from a
  distance between a normal random point and a tangent cone.  Different
  types of tangent cones will then be shown to lead to results (i) and
  (ii).
  
  Let $W(D_{r_n})$ be the set of points
  $(\bar\gamma,\bar\delta)\in(0,\infty)\times\mathbb{R}$ such that
  $D_{r_n}(\bar\gamma,\bar\delta)\le 0$.  Let
  \[
  \lambda_n = n\cdot \min_{(\bar\gamma,\bar\delta)\in W(D_{r_n})}
  (\hat\gamma-\bar\gamma)^2+(\hat\delta - \bar\delta)^2
  \]
  be the squared and scaled distance between the random point
  $(\hat\gamma,\hat\delta)$ and $W(D_{r_n})$.  In Figure \ref{fig:curves}
  the set $W(D_{r_n})$ corresponds to the region between and including the
  two curve branches.  The Behrens-Fisher likelihood equations have three
  distinct real solutions if and only if $D_{r_n}(\hat\gamma,\hat\delta)>0$
  if and only if $\lambda_n>0$.
  
  By the central limit theorem and the delta method, the two random
  variables $A_n=\sqrt{n}(\hat\gamma-\gamma)$ and
  $B_n=\sqrt{n}(\hat\delta-\bar\delta)$ converge jointly to a centered
  bivariate normal distribution $N_2(0,\Sigma)$.  In order to make use of
  this convergence, we rewrite
  \begin{align*}
    \lambda_n &= \min_{(\bar\gamma,\bar\delta)\in W(D_{r_n})}
    \big[A_n-\sqrt{n}(\bar\gamma-\gamma)\big]^2+
    \big[B_n - \sqrt{n}(\bar\delta-\delta)\big]^2.    
  \end{align*}
  The limits, for $n\to\infty$, of convergent sequences of the form
  $\sqrt{n}[ (\bar\gamma_n, \bar\delta_n)- (\gamma,\delta) ]$ with
  $(\bar\gamma_n,\bar\delta_n)\in W(D_{r_n})$ form the tangent cone
  $T(\gamma,\delta)$ of the set $W(D_r)$ at $(\gamma,\delta)$.  It thus
  follows from \citet[Lemma 7.13]{vandervaart} that as $n$ tends to
  infinity, the random distance $\lambda_n$ converges in distribution to
  the distance
  \begin{align*}
    \lambda_\infty = \min_{(\bar\gamma,\bar\delta)\in T(\gamma,\delta)}
    (Z_1-\bar\gamma)^2+ (Z_2-\bar d)^2,
  \end{align*}
  between the normal random vector $Z=(Z_1,Z_2)\sim N_2(0,\Sigma)$ and
  $T(\gamma,\delta)$.
  
  {\em Case (i):} If $(\gamma,\delta)$ is a non-singular point of $V(D_r)$,
  then $T(\gamma,\delta)$ is a half-space $H$ comprising all points on and
  to one side of a line through the origin.  The normal vector of this line
  is given by the gradient $\nabla D_r(\gamma,\delta)$.  The probability
  $P(\lambda_\infty>0)=P(\Sigma^{-1/2}Z\in\Sigma^{-1/2}H)$ is equal to
  $1/2$ because $\Sigma^{-1/2}Z\sim N_2(0,I)$ is standard normal and
  because $\Sigma^{-1/2}H$ is still a half-space with the origin on its
  boundary.  Since $P(D_{r_n}(\hat\gamma,\hat\delta)>0)=P(\lambda_n>0)$
  converges to $P(\lambda_\infty>0)$ we have established claim (i).
  
  {\em Case (ii):} If $(\gamma,\delta)$ is a singular point, then
  $T(\gamma,\delta)$ is all of $\mathbb{R}^2$; compare Figure
  \ref{fig:curves}.  Thus $P(\lambda_\infty>0)=0$, which implies claim
  (ii). 
\end{proof}

For the curious reader, we remark that the tangent cone to the curve
$V(D_r)$ at its singular point $(\gamma,\delta)$ with $\gamma,\delta>0$ is
the half-ray of points $(\gamma,\delta)$ with $\delta\ge 0$ and
$\gamma\sqrt{3(2r+1)}=\delta(r-1)$.  If $r=1$, then this half-ray is 
the non-negative $\delta$-axis.  
The half-ray has positive slope
if $r>1$.  The slope is negative if $r< 1$.

\begin{figure}[t]
  \centering
  \includegraphics[width=10cm]{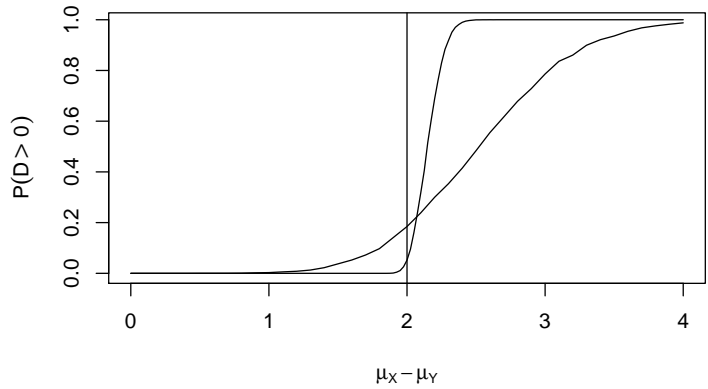}
  \caption{Simulations for the probability of three distinct real roots to the
    Behrens-Fisher likelihood equations ($\sigma_X^2=\sigma_Y^2=1$). The
    steeper curve corresponds to $n=m=1000$; the other to $n=m=15$.}
  \label{fig:asy-alternative}
\end{figure}

We conclude by illustrating the results obtained in this section in Figure
\ref{fig:asy-alternative}, which shows simulations on the probability of
three distinct real roots to the Behrens-Fisher likelihood equations.  This
figure addresses the case $\gamma=r=1$ in which $\delta=\mu_X-\mu_Y$.  The
simulations 
confirm Theorem \ref{thm:critical-smooth}(ii) because the probability of
three distinct real roots appears to converge to zero if $\delta=2$ and
$(\gamma,\delta)=(1,2)$ is a singularity of $V(D_1)$.

\bibliographystyle{plainnat}

\begin{thebibliography}{1}
\providecommand{\natexlab}[1]{#1}
\providecommand{\url}[1]{\texttt{#1}}
\expandafter\ifx\csname urlstyle\endcsname\relax
  \providecommand{\doi}[1]{doi: #1}\else
  \providecommand{\doi}{doi: \begingroup \urlstyle{rm}\Url}\fi

\bibitem[Basu et~al.(2003)]{basu:2003} {Basu, S., Pollack, R., Roy,
  M.-F.}, 2003. {\em Algorithms in Real Algebraic Geometry}, Berlin:
  Springer-Verlag. 
  





\bibitem[Jensen(1992)]{jensen:1992} Jensen, J. L., 1992.
The modified signed likelihood statistic and saddlepoint
approximations. {\em Biometrika} 79, 693--703. 


\bibitem[Pachter and  Sturmfels(2005)]{ascb} { Pachter, L., and
    Sturmfels, B.}, 2005. {\em Algebraic Statistics for Computational Biology}, 
  Cambridge: Cambridge University Press.
  
\bibitem[Sugiura and  Gupta(1987)]{sugiura} { Sugiura, N., and Gupta,
    A. K.}, 1987. Maximum likelihood estimates for Behrens-Fisher problem, {\em
    Journal of the Japan Statistical Society}, 17, 55--60.
  
\bibitem[van der Vaart(1998)]{vandervaart} { van der Vaart, A. W.}, 1998. {\em
    Asymptotic Statistics}, 
  Cambridge: Cambridge University Press.



\end{thebibliography}



\bigskip

\end{document}